\documentclass[12pt]{amsart}
\usepackage{amsmath}
\usepackage{amsfonts}
\usepackage{amssymb}

\newcommand{\be}{\begin{eqnarray}}
\newcommand{\ee}{\end{eqnarray}}
\newcommand{\supp}{\mbox{\rm supp}}
\newcommand{\muhat}{{\widehat{\mu}}}
\newcommand{\card}{\#}
\newcommand{\const}{\mbox{\rm const}}

\newcommand{\half}{\frac{1}{2}}

\newcommand{\de}{\delta}
\newcommand{\e}{{\varepsilon}}
\newcommand{\Th}{{\theta}}

\newcommand{\R}{{\mathbb R}}
\newcommand{\Q}{{\mathbb Q}}
\newcommand{\Z}{{\mathbb Z}}
\newcommand{\Nat}{{\mathbb N}}

\newcommand{\Fk}{{\mathcal F}}

\newcommand{\jk}{J^{(k)}}

\newcommand{\La}{{\Lambda}}

\newcommand{\om}{\omega}

\renewcommand{\th}{\theta}
\newcommand{\ga}{\gamma}
\newcommand{\Om}{\Omega}

\newcommand{\BR}{\mathbb{R}}
\newcommand{\BZ}{\mathbb{Z}}
\newcommand{\BQ}{\mathbb{Q}}

\newcommand{\wh}{\widehat}

\renewcommand{\P}{\mathcal{P}}

\renewcommand{\sp}{\mbox{sp}}

\renewcommand{\phi}{\varphi}

\def\smfrac#1/#2{\leavevmode\kern-.1em
    \raise.5ex\hbox{\the\scriptfont0 #1}\kern-.1em
    /\kern-.15em\lower.25ex\hbox{\the\scriptfont0 #2}}

\newtheorem{lemma}{Lemma}[section]

\newtheorem{prop}[lemma]{Proposition}
\newtheorem{thm}[lemma]{Theorem}
\newtheorem{cor}[lemma]{Corollary}
\theoremstyle{definition}

\theoremstyle{remark}
\newtheorem{rmk}[lemma]{Remark}

\numberwithin{equation}{section}

\begin{document}

\title[Spectra of Bernoulli convolutions as multipliers]
{Spectra of Bernoulli convolutions as multipliers in $L^p$ on
the circle }
\author{Nikita Sidorov}
\address{Department of Mathematics, UMIST, P.O. Box 88,
Manchester M60 1QD, United Kingdom. E-mail:
Nikita.A.Sidorov@umist.ac.uk}
\author{Boris Solomyak}
\address{Box 354350, Department of Mathematics, University of
Washington, Seattle, WA~98195, USA. E-mail:
solomyak@math.washington.edu}
\date{\today}

\thanks{The first author was supported by the EPSRC grant no
GR/R61451/01. The second author was supported in part by the NSF
Grant DMS~0099814.}



\begin{abstract} It is shown that the closure of the set of
Fourier coefficients of the Bernoulli convolution $\mu_\theta$
parameterized by a Pisot number $\theta$, is countable. Combined
with results of Salem and Sarnak, this proves that for every
fixed $\theta>1$ the spectrum of the convolution operator
$f\mapsto \mu_\theta*f$ in $L^p(S^1)$ (where $S^1$ is the circle
group) is countable and is the same for all $p\in(1,\infty)$, namely,
$\overline{\{\widehat{\mu_\theta}(n) : n\in\mathbb{Z}\}}$.
Our result answers the question raised by P.~Sarnak in \cite{Sar}.
We also consider the sets $\overline{\{\widehat{\mu_\theta}(rn)
: n\in\mathbb{Z}\}}$ for $r>0$ which correspond to a linear
change of variable for the measure. We show that such a set is
still countable for all $r\in\Q(\theta)$ but uncountable (a
non-empty interval) for Lebesgue-a.e.\ $r>0$.
\end{abstract}

\maketitle

\section{Introduction and main results}\label{intro}

Let $\nu$ be a Borel probability measure on $S^1=\BR/\BZ$, and let
\[
T_\nu:L^p(S^1)\to L^p(S^1)
\]
be the convolution operator, namely, $T_\nu f:=\nu*f$. Put
\[
\Fk_\nu:=\{\wh\nu(n) : n\in\BZ\}\,.
\]
It is shown by Sarnak \cite{Sar} that if the closure
$\overline{\Fk_\nu}$ has capacity zero, then the following
identity relation for the spectra of $T_\nu$ in different
$L^p(S^1)$ is satisfied:
\begin{equation}
\sp(T_\nu, L^p)\equiv \sp(T_\nu, L^2)=\overline{\Fk_\nu},\quad
1<p<+\infty. \label{spec}
\end{equation}

The present paper deals with the case when $\nu$ is a Bernoulli
convolution. Recall that for any $\th>1$ the {\em Bernoulli
convolution} parameterized by $\th$ is defined as follows:
\[
\mu_\th:=*\prod_{k=0}^\infty\left(\frac12\,\de_{-\th^{-k}}+
\frac12\,\de_{\th^{-k}}\right)
\]
(where $\de_x$ denotes the Dirac $\de$-measure at $x$). Thus,
\begin{equation} \label{eq-four}
\wh{\mu_\th}(t) = \prod_{k=0}^\infty \cos(2\pi \th^{-k}t).
\end{equation}
We can view $\mu_\th$ as a measure on the line, i.e., $t\in \R$ in
(\ref{eq-four}). The induced measure on the circle has Fourier
coefficients $\{\wh{\mu_\th}(n) : n\in \Z\}$.

As is well known, $\supp\,\mu_\th\subset
I_\th=[-\frac{\th}{\th-1},\frac{\th}{\th-1}]$ for any $\th>1$.
Moreover, for $\th>2$ the measure $\mu_\th$ is usually called the
{\em Cantor-Lebesgue measure} (parameterized by $\th$), and its
support is the Cantor set with constant dissection ratio $\th$. On
the other hand, $\supp\,\mu_\th=I_\th$ for $\th\in (1,2]$.

In \cite[Th~II, p.~40]{salem1} it is shown that unless $\th$ is
a {\em Pisot number} (an algebraic integer greater than 1 whose
conjugates are all less than 1 in modulus),
$\wh{\mu_\th}(t)\to0$ as $t\to+\infty$ along the reals, whence
$\overline{\Fk_{\mu_\th}}$ is countable and thus, (\ref{spec})
is satisfied.

Sarnak \cite{Sar} considered the case of the classical
Cantor-Lebesgue measure ($\th=3$), for which he proved that
$\overline{\Fk_{\mu_\th}}$ is countable and therefore (\ref{spec})
holds. As stated in \cite{Sar}, the same approach can be
applied to the case of an arbitrary integer $\th\ge3$, and the
only case left is the irrational Pisot numbers $\th$. The question
about the limit points of the Fourier coefficients for this class
of measures was raised in \cite{Sar}.


\begin{thm}\label{main} The set of limit points of the sequence
$\{\wh{\mu_\th}(n): n\in\BZ\}$, with an irrational Pisot parameter
$\th$, is countable, so (\ref{spec}) holds for $\nu = \mu_\th$.
\end{thm}

\begin{rmk}\label{lp-imp} There is apparently another way to
obtain (\ref{spec}) for $\nu = \mu_\th$, without getting
countability of the spectrum. It is known that if a Borel
probability measure $\nu$ on $S^1$ is $L^p$-improving\footnote{A
measure $\nu$ is called {\em $L^p$-improving} for some
$p\in(1,\infty)$, if there exists $q=q(p)>p$ such that $\nu*f\in
L^q(S^1)$ for any $f\in L^p(S^1)$. If $\nu$ improves some $L^p$,
then it improves all of them for $1<p<+\infty$, see
\cite{GHR}.}, then, like in the case in question, (\ref{spec})
holds, see \cite[Th. 4.1]{GHR}. Christ \cite{Christ} proved that
$\mu_\theta$ is $L^p$-improving for all $\th>2$, and he made a
remark that the same argument works for $\th\in (1,2)$ as well.
\end{rmk}

Since we are considering measures on the circle, one may argue
that it is in fact more natural to have the Bernoulli convolution
measure supported on an interval of length~1, rather than on
$I_\th$ (whose length is $2\th/(\th-1)>1$). This is achieved by a
linear change of variable, resulting in the Fourier coefficients
$\{\muhat((\frac{\th-1}{2\th})n):\, n \in \Z\}$, and we show that
the analog of Theorem~\ref{main} is still valid (see below).

More generally, one may inquire what happens under an arbitrary
scale change. It turns out that the situation is rather delicate.

For $r>0$ put
$$
\Fk_{\nu,\,r} := \{\wh\nu(rn) : n\in\BZ\}.
$$
For a set $E\subset\R$ let $E'$ denote the derived set of $E$,
that is, the set of its limit points.

\begin{thm}\label{thm2} Let $\th\neq2$ be a Pisot number. Then
\begin{enumerate}
\item for any positive $r\in \Q(\th)$, the set
$\Fk_{\mu_\th,r}'$ is countable;
\item for Lebesgue a.e.\ $r>0$, the set
$\Fk_{\mu_\th,r}'$ is a non-empty interval.
\end{enumerate}
\end{thm}

\begin{cor} For any $\th>2$, the spectrum in $L^p(S^1),\ p>1$,
of the convolution operator corresponding to the Cantor-Lebesgue
measure with the constant dissection ratio $\th$, constructed on
$[0,1) \simeq S^1 = \R/\Z$, is countable.
\end{cor}

This follows from Theorem~\ref{thm2}~(i), since
$\frac{\th-1}{2\th}\in \Q(\th)$ and translating the measure by
$\half$ results in multiplying the Fourier coefficients by
$(-1)^n$.


Denote $E^{(1)}:=E',\ E^{(n+1)}:=(E^{(n)})'$. For $r\in \Q(\th)$,
it is natural to ask what is the cardinality of the second, third,
etc., derived sets for $\Fk_{\mu_\th,r}$. The following theorem
answers this question.

\begin{thm}\label{thm:deriv}
For any Pisot $\th\neq2$ and any positive $r\in\Q(\th)$ the set
$\Fk_{\mu_\th,r}^{(n)}$ is countable for each $n\ge1$.
\end{thm}

\section{Proof of countability: the model case}\label{sec2}

The core of the paper is the proof of Theorem~\ref{thm2}~(i)
(which, of course, implies Theorem~\ref{main}). Our proof is
loosely based on the method used in \cite{Sar} for $\th=3$. On
the other hand, the case of irrational $\th$ requires extra
tools more common for the theory of Pisot numbers (in the spirit
of the monograph \cite[Chapter~VIII]{Cas}).

The structure of the rest of the paper is as follows: in this and
the next sections we are going to show that
\begin{equation}
\mbox{card}\,\Fk_{\mu_\th,r}'\le\aleph_0,\quad
r\in\Q(\th)\cap(0,\infty).\label{eq:loeq}
\end{equation}
Note that our proof applies to the case $\th\in\Nat$
as well. Combined with Theorem~\ref{thm:deriv}
(proved in Section~\ref{sec:proofs}) this yields
Theorem~\ref{thm2}~(i). Theorem~\ref{thm2}~(ii) is also proved in
Section~\ref{sec:proofs}.

The Pisot numbers $\th\ne 2$ are distinguished by the fact that
$\muhat_\th(t) \not\to 0$, as $t\to +\infty$ over the reals.
This was proved by Erd\H{o}s \cite{E} (for integers $\th>2$ this
had been known earlier). It is easy to see that also
$\muhat_\th(n) \not \to 0$, as $n\to +\infty$ over the integers,
see Section 4. If $\th=2$, then $\mu_\th$ is absolutely
continuous, so there is nothing to prove.

For the rest of the paper, we fix a Pisot number $\th\ne 2$ and denote by
$\th_2,\dots,\th_m$ the conjugates of $\th=\th_1$. Since $\th$ is Pisot,
\begin{equation}
\rho=\max_{i\ge2}\,|\th_i|\in[0,1)\label{defro}
\end{equation}
($\rho=0$ if and only if $\th\in\Nat$).

Let us begin the proof of the inequality~(\ref{eq:loeq}).
Denote by $\langle\cdot\rangle$, $\|\cdot\|$ the nearest integer and the
distance to the nearest integer respectively.



To simplify notation, denote $\mu:=\mu_\th$. It suffices
to prove that there are at most countable many limit points for
the set $\{|\muhat(rn)|:\ n\in \Nat\}$. Fix $\eta>0$ and assume
that integers $n_k \to +\infty$ are such that
\begin{equation} \label{eq-A}
|\muhat(rn_k)| \to a, \quad a \ge \eta.
\end{equation}
Similarly to \cite{Sar}, our goal is to show that there can be
only a countable set of such $a$'s for any fixed $\eta$; this will
yield (\ref{eq:loeq}).

There exist $N_k \in \Z$ and $y_k \in [1,\th)$ such that
\begin{equation} \label{def-yk}
y_k = 2rn_k \theta^{-N_k}.
\end{equation}
Let
\begin{equation}\label{eq:Kj}
y_k \Th^j = K_j^{(k)} + \de_j^{(k)},\quad j = 1,\ldots,N_k,
\end{equation}
where $\de_j^{(k)} \in (-\half,\half]$ and $K_j^{(k)}=\langle
y_k\Th^j\rangle\in \Nat$. By (\ref{eq-four}) and (\ref{def-yk}),
\begin{equation}\label{eq-B}
|\muhat(rn_k)| = \prod_{j=-\infty}^{N_k} |\cos(\pi y_k \th^j)| =
\prod_{j=-\infty}^{N_k} \cos\bigl(\pi \de_j^{(k)}\bigr).
\end{equation}
Let $m$ be the degree of $\th$, with the minimal polynomial
$x^m-d_1x^{m-1}-\dots-d_m$. For the rest of the paper we fix a
$\de$ which satisfies
\begin{equation}
0<\delta < (1+|d_1| + \cdots + |d_m|)^{-1}. \label{delta}
\end{equation}
The reason for the choice of $\de$ is the following

\begin{lemma} \label{lem-recur}
Suppose that $|\de_j^{(k)}|\le \de$ for $j = A_k + 1, \ldots, A_k
+ b$, where $0 \le A_k \le N_k -b$ and $b > m$. Then \be
\label{eq-recur} K_{j+m}^{(k)} = d_1 K_{j+m-1}^{(k)} + \cdots +
d_m K_{j}^{(k)}, \ee for $j = A_k +1, \ldots, A_k + b-m$.
\end{lemma}
\begin{proof} By our condition and (\ref{delta}), for $j\ge
A_k+1$,
\[
|K_{j+m}^{(k)}-d_1 K_{j+m-1}^{(k)}-\dots-d_m
K_{j}^{(k)}|\le\de(1+|d_1|+\dots+|d_m|) <1.
\]
As $K_j^{(k)}$'s and $d_i$'s are integers, we are done.
\end{proof}

We want to estimate the number of $\de_j^{(k)}$'s that are greater
than $\de$ in modulus. Let $L \in \Nat$ be such that $(\cos(\pi
\delta))^L \le \eta/2$. It follows from (\ref{eq-A}) and
(\ref{eq-B}) that for $k$ sufficiently large,
\begin{equation} \label{def-Lk}
L_k:=\card \bigl\{j \in [1, N_k]:\ \bigl|\de_j^{(k)}\bigr|> \de
\bigr\} \le L.
\end{equation}
Since we only care about the limit, we can assume without loss of
generality that (\ref{def-Lk}) holds for all $k$.

The rest of the proof is somewhat technical, so we believe that it
is helpful first to present a sketch in the special model case
$L_k=1$ and $r=1$. This will be done in the rest of the section.

Thus, let us assume for the moment that
$\bigl|\de_j^{(k)}\bigr|\le \de$ for all $j = 1,\ldots,N_k$,
except possibly $j= J_k$. There are three possibilities: (a) $\sup
J_k< \infty$, (b) $\sup(N_k- J_k) < \infty$, and (c) $\sup J_k
=\infty$ and $\sup(N_k - J_k) = \infty$. By passing to a
subsequence, we can assume that we actually have one of the
following cases:

\smallskip

\textbf{Case 1:} $J_k = R$ (independent of $k$);

\smallskip

\textbf{Case 2:} $J_k = N_k - R$;

\smallskip

\textbf{Case 3:} $J_k \to \infty$ and $N_k - J_k \to \infty$.

\medskip
\noindent \textbf{Case 1.} By Lemma~\ref{lem-recur}, the sequence
$\{K_j^{(k)}\}$ satisfies the recurrence relation (\ref{eq-recur})
for $j = R+1,\ldots, N_k-m$ (for $k$ large enough to satisfy
$N_k>R+m$). Then we can express $K_j^{(k)}$ in terms of $\theta$
and its conjugates $\Th_2,\ldots,\Th_m$ as follows:
\begin{equation} \label{eq-exp1}
K_j^{(k)} = c_1^{(k)} \Th^j + \sum_{i=2}^m c_i^{(k)} \Th_i^j,\quad
j = R+1,\ldots,N_k.
\end{equation}
Observe that the coefficients $c_i^{(k)}$ are completely
determined by $K_j^{(k)}$ for $j = R+1,\ldots,R+m$. These
$K_j^{(k)}$'s are integers bounded by $\Th^{R+m+1} + 1$ (as $y_k
\le \Th$ and $K_j^{(k)}$ is the nearest integer to $y_k \Th^j$).
Thus, there are finitely many possibilities for $c_i^{(k)}$ and we
can assume, passing to a subsequence, that $c_i^{(k)}=c_i$ do not
depend on $k$.

Let $y:=c_1$; the first important point is that $y\in \Q(\th)$.
This follows from the Cramer's Rule, solving the linear system
(\ref{eq-exp1}), with $j=R+1,\ldots,R+m$, for $c_i$.
Alternatively, note that $\|y\th^j\|\to 0$ as $j\to+\infty$ by
(\ref{eq-exp1}), and the fact that $y\in \Q(\th)$ is a part of the
well-known Pisot-Vijayaraghavan theorem (see \cite{Cas}).

Now comes the crucial point---we have to use that $n_k$ is an
integer\footnote{Note that the set of limit points of $\muhat(t)$,
as $t\to\infty$ over the reals, is an interval---see
Lemma~\ref{lem:int} below.}. We assumed that $r=1$, so $n_k =
\half y_k \th^{N_k}$ from (\ref{def-yk}). We have $2n_k =
K_{N_k}^{(k)}$, since both sides are integers. Hence by
(\ref{eq:Kj}) and (\ref{eq-exp1}),
$$
y_k \Th^{N_k} = K_{N_k}^{(k)}=y\Th^{N_k} + O(\rho^{N_k}),
$$
where $\rho$ is given by (\ref{defro}) and the implied constant in
$O$ is independent of $k$. Thus, $y_k=y+O(\Th^{-N_k}\rho^{N_k})$,
and an elementary argument yields
$$
|\muhat(n_k)| \to \prod_{j=-\infty}^\infty |\cos(\pi y \Th^j)|.
$$
(A more general statement is proved below, in
Lemma~\ref{lem-converge}.) Since the right-hand side depends only
on $y \in \Q(\Th)$, the number of possible limit points in this
case is at most countable.

\medskip\noindent
\textbf{Case 2.} By Lemma~\ref{lem-recur}, $K_j^{(k)}$'s satisfy
the recurrence relation (\ref{eq-recur}) for $j = 1,\ldots, N_k -
R-m-1$, when $k$ is sufficiently large. Passing to a subsequence,
we can assume that $K_j^{(k)}=K_j$ do not depend on $k$ for $j \le
N_k - R-1$ and
\begin{equation}\label{eq-Kjj}
K_j = y \Th^j + \sum_{i=2}^m c_i \Th_i^j,\quad j = 1,\ldots, N_k -
R-1.
\end{equation}
Again we have $y \in \Q(\th)$. Extend $K_j$ by (\ref{eq-Kjj}) to
$j = N_k - R,\ldots,N_k$; in other words, we extend $K_j$ to
satisfy the recurrence relation (\ref{eq-recur}). We cannot claim
that $K_j = K_j^{(k)}$ for $j = N_k - R,\ldots,N_k$; however, it
is easy to see from the recurrence that
$$
\bigl|K_j^{(k)}-K_j\bigr| \le C_R,\quad j = N_k - R,\ldots,N_k,
$$
where $C_R$ does not depend on $k$. (This is proved below, in
Lemma~\ref{lem-recur2}.) Again, passing to a subsequence, we can
assume that $K_{N_k}^{(k)}-K_{N_k} = A$ is a constant. Using that
$n_k$ is an integer, we obtain
$$
y_k \Th^{N_k} = K_{N_k}^{(k)} = y \Th^{N_K} + A + O(\rho^{N_k}),
$$
where the implied constant in $O$ is independent of $k$. Thus $y_k
= y + A\Th^{-N_k} + O(\Th^{-N_k}\rho^{N_k})$, and it is not hard
to show that
$$
|\muhat(n_k)| \to \prod_{j=-\infty}^\infty |\cos(\pi y
\Th^j)|\cdot \prod_{j=0}^\infty |\cos(\pi A \Th^{-j})|.
$$
(A more general statement will be proved below, in
Lemma~\ref{lem-converge}.) Since $y \in \Q(\Th)$ and $A\in \Z$,
the number of possible limit points in this case is again at most
countable.

\medskip
\noindent \textbf{Case 3.}  By Lemma~\ref{lem-recur},
$\bigl\{K_j^{(k)}\bigr\}_{j=1}^{J_k-m-1}$ and
$\bigl\{K_j^{(k)}\bigr\}_{j=J_k+1}^{N_k-m}$ satisfy the recurrence
relation (\ref{eq-recur}). As in Case~2, we can assume by passing
to a subsequence that $K_j^{(k)}=K_j$ for $j=1,\dots,J_k-1$,
whence
$$
K_j = y \Th^j + \sum_{i=2}^m c_i \Th_i^j,\quad j = 1,\ldots,
J_k-1,
$$
for some $y\in \Q(\Th), c_2,\ldots, c_m$. Also, as in Case~2, we
extend $K_j$ to $j \ge J_k$ to satisfy the same recurrence
relation and check that
$$
\bigl|K_j^{(k)}-K_j\bigr| \le C_R,\quad j =J_k+1,\ldots,J_k+m.
$$
Let $S_j^{(k)} = K_{J_k+j}^{(k)}-K_{J_k+j}$ for
$j=1,\ldots,N_k-J_k$. Hence there exist $b_i^{(k)}\in\BQ(\th_i)$,
with $i=1,\ldots,m$, such that
$$
S_j^{(k)} = b_1^{(k)}\Th^j + \sum_{i=2}^m b_i^{(k)} \Th_i^j,\quad
j = 1,\ldots,m.
$$
Because of the bounds on $S_j^{(k)}$, there are finitely many
possibilities for $b_i^{(k)}$, so we can assume that they do not
depend on $k$, passing to a subsequence. Let $z = b_1 =
b_1^{(k)}$. We have $S_j = S_j^{(k)}= z\Th^j + O(\rho^j)$ for $j
\ge 1$. Now observe that
$$
K_j^{(k)} = K_j + S_{j-J_k},\ \ \ j=J_k+1,\ldots,N_k.
$$
Thus, using that $n_k$ is an integer, we obtain
$$
y_k \Th^{N_k} = K_{N_k}^{(k)} = K_{N_k} + S_{N_k-J_k} = y\Th^{N_k}
+ z\Th^{N_k-J_k} + O(\rho^{N_k-J_k}),
$$
where the implied constant in $O$ is independent of $k$. Since
$J_k\to \infty$ and $N_k-J_k\to \infty$,  is not hard to show that
$$
|\muhat(n_k)| \to \prod_{j=-\infty}^\infty |\cos(\pi y
\Th^j)|\cdot \prod_{j=-\infty}^\infty |\cos(\pi z \Th^j)|.
$$
(We will prove a more general statement below, in
Lemma~\ref{lem-converge}.) As $y,z\in \Q(\Th)$, the number of
possible limit points in this case is at most countable.

This concludes the sketch of the proof of (\ref{eq:loeq}) in the
model case $L_k=1$ and $r=1$. The idea for the general case is as
follows: we gather all indices $j$, for which
$\bigl|\de_j^{(k)}\bigr|>\de$, in groups in such a way that the
distance between any two adjacent groups goes to the infinity as
$k\to\infty$. Then we treat each group similarly to one of the
three cases considered in this section, depending on the position
of this group (``beginning", ``middle" or ``end") and finally,
justify passing to the limit in the key Lemma~\ref{lem-converge}.

\section{Proof of countability: the general case}

We continue with the proof of the general case where we left it,
after the definition of $L_k$ (\ref{def-Lk}). Let $1\le I_1^{(k)}
< I_2^{(k)} < \ldots < I_{L_k}^{(k)} \le N_k$ be all the indices
$j$ for which $\bigl|\de_j^{(k)}\bigr|>\de$. Since $L_k \le L$, we
can assume that $L_k = L'$ does not depend on $k$, passing to a
subsequence. Further, passing to a subsequence, we can assume that
for all $i=1,\ldots, L'-1$, either $I_{i+1}^{(k)}-I_i^{(k)} = R_i$
(independent of $k$), or $I_{i+1}^{(k)}-I_i^{(k)} \to \infty$, as
$k\to \infty$. Also, either $I_1^{(k)} = R_0$ or
$I_1^{(k)}\to\infty$ and either $N_k - I_{L'}^{(k)} = R_{L'}$ or
$N_k - I_{L'}^{(k)}\to\infty$. Let
$$
R = \max\{R_i:\ i= 0,\ldots,L'\} + 1.
$$
We can find $M\in\{1,\dots,L'\}$ and integers
$$
1 = J_0^{(k)} < J_1^{(k)} < \ldots < J_M^{(k)} < J_{M+1}^{(k)} =
N_k
$$
so that $J_{i+1}^{(k)} - J_i^{(k)} \to\infty$ for $i=0,\dots,M$,
and
$$
\bigl|\de_j^{(k)}\bigr| \le \de\quad\mbox{for all}\ j \in
\{1,\ldots,N_k\}\quad\mbox{such that}\ \min_i
\bigl|j-J_i^{(k)}\bigr| \ge R.
$$
By Lemma~\ref{lem-recur}, $\bigl\{K_j^{(k)}\bigr\}$ satisfy the
recurrence relation (\ref{eq-recur}) for $J_i^{(k)} + R \le j \le
J_{i+1}^{(k)} - R-m$, with $i = 0,\ldots, M$ (for $k$ large enough
to satisfy $J_{i+1}^{(k)}-J_i^{(k)}>2R+m$). In particular, this is
true for $1+R \le j \le \jk_1-R-m$. Thus we can write
$$
K_j^{(k)} = c_1^{(k)} \Th^j + \sum_{i=2}^m c_i^{(k)} \Th_i^j,\quad
j = R+1,\ldots,\jk_1-R
$$
(for $k$ sufficiently large). The coefficients $c_i^{(k)}$ are
completely determined by $K_j^{(k)}$ for $j = R+1,\ldots,R+m$,
which are integers bounded by $\Th^{R+m+1} + 1$. Thus, there are
finitely many possibilities for $c_i^{(k)},\ i=1,\dots,m$. Hence
we can assume, passing to a subsequence, that $c_i^{(k)}=c_i$ do
not depend on $k$. Thus, $K_j^{(k)}=K_j$ for $j\le \jk_1-R$.
Denote $z_0:=c_1$, so that
$$
K_j = z_0 \Th^j + O(\rho^j).
$$
As in Case~1, we have $z_0 \in \Q(\th)$ ($z_0$ is a natural analog
of $y$ from the previous section). Next we repeat the argument
from Case~3. Extend the sequence $\{K_j\}$ to $j > \jk_1-R$ so
that it satisfies the recurrence relation for all $j$. We need the
following lemma.

\begin{lemma} \label{lem-recur2}
Suppose that $K_j^{(k)} = K_j$ for $j \le A_k$ and $\{K_j\}$
satisfies the recurrence relation (\ref{eq-recur}) for all $j$.
Then for any $p\in\Nat$ there exists $C_p=C_p(\th)>0$ such that
\be \label{eq-est} |K_j^{(k)}-K_j| \le C_p,\quad j = A_k+1,\ldots,
A_k+p. \ee
\end{lemma}
\begin{proof}
This is proved by induction. Since $K_j^{(k)}$ is the nearest
integer to $y_k \Th^j$, for any $j\ge 1$,
$$
|K_{j+m}^{(k)} - d_1 K_{j+m-1}^{(k)} - \dots - d_m
K_j^{(k)}|\le\frac12(1+ |d_1| + \cdots + |d_m|) =:\gamma.
$$
It follows that we can take $C_1=\ga$ in (\ref{eq-est}). Suppose
(\ref{eq-est}) is verified for some $p$. Then we have for
$j=A_k+p+1$,
\begin{align*}
|K_{j+1}^{(k)}-K_{j+1}|&=|K_{j+1}^{(k)}-d_1 K_j^{(k)}-\dots-d_m
K_{j-m+1}^{(k)}\\
&+d_1 K_j^{(k)}+\dots+d_m K_{j-m+1}^{(k)}-d_1K_j-\dots-d_m
K_{j-m+1}|\\
&\le \ga+(|d_1|+\dots+|d_m|)C_p < \ga(1+2C_p).
\end{align*}
Thus, we may put $C_{p+1} = \ga(1+2C_p)$, and the lemma is proved.
\end{proof}
Let
$$
S_j^{(k)} = K^{(k)}_{\jk_1+j}-K_{\jk_1+j},\quad j=R+1,\ldots,R+m.
$$
By Lemma~\ref{lem-recur2},
$$
\bigl|S_j^{(k)}\bigr| \le C_{2R+m},\quad j = R+1, \ldots, R+m.
$$
Therefore, we can assume (passing to a subsequence) that
$S_j^{(k)}$'s do not depend on $k$ for $j=R+1,\ldots,R+m$. We can
find $z_1 \in \Q(\Th), c_2',\ldots,c_m'$ so that \be
\label{eq-exp2} S_j = S_j^{(k)} = z_1\Th^j + \sum_{i=2}^m c_i'
\Th_i^j,\ \ \ j = R+1,\ldots,R+m. \ee Extend $S_j$ to $j>R+m$ by
the formula (\ref{eq-exp2}), so that they satisfy the recurrence
relation. Now observe that
$$
K_j^{(k)} = K_j + S_{j - \jk_1},\quad j=\jk_1+R+1, \ldots, \jk_2-R
$$
(for $k$ large enough to satisfy $J_2^{(k)}-J_1^{(k)}>2R+m$), as
both sides agree for $j= \jk_1+R+1, \ldots, \jk_1+R+m$, and
satisfy the same recurrence relation of length $m$. It follows
that
$$
K_j^{(k)} = z_0 \Th^j + z_1 \Th^{j-\jk_1} +
O\bigl(\rho^{j-\jk_1}\bigr),\quad j= \jk_1+R+1,\ldots, \jk_2-R,
$$
where the implied constant in $O$ is independent of $k$.

Next we repeat the same argument and obtain, by induction, that
for $i=2,\ldots,M$,
\begin{equation}
\begin{aligned}
K_j^{(k)} = &z_0 \Th^j + z_1 \Th^{j-\jk_1} + \cdots + z_i
\Th^{j-\jk_i} + O(\rho^{j-\jk_i}),\\
&j = \jk_i+R,\ldots,\jk_{i+1}-R. \label{eq-exp3}
\end{aligned}
\end{equation}
Indeed, for each $i$ extend $K_j^{(k)}$ from $j< \jk_i-R$ to
larger $j$'s by recurrence; denote them $Q_j^{(k)}$. Put
$$
T_j^{(k)} := K^{(k)}_{\jk_i+j}-Q^{(k)}_{\jk_i+j},\quad
j=R+1,\ldots,R+m.
$$
Then $\bigl|T_j^{(k)}\bigr|\le C_{2R+m}$ for $j = R+1,\ldots,R+m$.
We can write $T_j^{(k)}$ as a linear combination of $\Th^j$ and
the powers of its conjugates.  As above, there are finitely many
possibilities for the coefficients (as $k$ varies), so we can
assume without loss of generality that they do not depend on $k$.
The coefficient at $\Th^j$ will be denoted by $z_i$, which yields
(\ref{eq-exp3}).

For $i=M$ the formula (\ref{eq-exp3}) becomes
\begin{equation}
\label{eq-KL}
\begin{aligned}
K_{N_k-R+j}^{(k)} = &z_0 \Th^{N_k-R+j} + \sum_{i=1}^M z_i
\Th^{N_k -\jk_i-R+j} + O\bigl(\rho^{N_k-\jk_M-R+j}\bigr),\\
&j=\jk_M+2R-N_k,\ldots,0
\end{aligned}
\end{equation}
(recall that $J_M^{(k)}-N_k \to -\infty$). As usual, the implied
constant in $O$ is independent of $k$. One last time extend
$K^{(k)}_j$ by recurrence, to $j=N_k-R+1,\ldots,N_k,\ldots$ Denote
the resulting integer sequence by
$\{L_{N_k-R+j}^{(k)}\}_{j=1}^\infty$. By Lemma~\ref{lem-recur2},
\begin{equation} \label{equ1}
|K_{N_k}^{(k)} - L_{N_k}^{(k)}| \le C_R.
\end{equation}
We have
\begin{equation}\label{eq-2rnk}
2rn_k = y_k \th^{N_k} = K_{N_k}^{(k)} + \delta_{N_k}^{(k)}.
\end{equation}
Now it's the time to use the fact that $n_k$ is an integer; this
is slightly more complicated than in Section~\ref{sec2}, because
the left-hand side of (\ref{eq-2rnk}) need not be an integer.
However, as $2r \in \Q(\th)\setminus\{0\}$ by assumption, we can
invert it in $\Q(\th)$, i.e.,
\begin{equation} \label{invert}
(2r)^{-1} = a_0 + a_1\th + \cdots + a_{m-1}\th^{m-1},
\end{equation}
for some $a_i \in \Q$. Thus,
\begin{equation}
\begin{aligned}\label{equ2}
n_k &=  (2r)^{-1}\bigl(K_{N_k}^{(k)}
+\delta_{N_k}^{(k)}\bigr)\\
    &= (a_0 + a_1\th + \cdots + a_{m-1}\th^{m-1}) K_{N_k}^{(k)}
+ (2r)^{-1} \delta_{N_k}^{(k)} \\
    & = a_0 L_{N_k}^{(k)} + a_1 L_{N_k+1}^{(k)} + \cdots +
a_{m-1} L_{N_K+m-1}^{(k)} + A_k.
\end{aligned}
\end{equation}
Let us estimate the ``error term'' $A_k$. By the definition of the
integers $L_{N_k-R+j}^{(k)}$, they satisfy (\ref{eq-KL}), with $K$
replaced by $L$, for $j=1,2,\ldots$ In particular,
\begin{equation}
\label{extra} L_{N_k+j}^{(k)} = z_0 \Th^{N_k+j} + \sum_{i=1}^M z_i
\Th^{N_k -\jk_i+j} + O(\rho^{N_k-\jk_M+j}),\quad j \ge 0.
\end{equation}
Hence in view of $N_k-\jk_M\to+\infty$,
$$
\lim_{k\to\infty}\bigl|L_{N_k+j+1}^{(k)} - \th
L_{N_k+j}^{(k)}\bigr| = 0, \quad 0\le j\le m-1,
$$
whence
$$
\left|L_{N_k}^{(k)}\sum_{j=0}^{m-1}a_{j}\th^{j} -
\sum_{j=0}^{m-1}a_j L_{N_k+j}^{(k)}\right|  \to 0,\quad
k\to+\infty.
$$
By (\ref{equ2}), (\ref{equ1}) and in view of
$\bigl|\de_{N_k}^{(k)}\bigr|\le\frac12$, we have $|A_k|\le C'$
for some constant $C'$ independent of $k$; namely, one may put
for $k$ large enough,
\[
C'= 1 + (4r)^{-1}(2C_R + 1),
\]
where $C_R$ is as in (\ref{equ1}).

On the other hand, it follows from (\ref{equ2}) that $A_k \in
s^{-1}\Z$ for some $s\in \Nat$ independent of $k$, because $a_i
\in \Q$ and the $n_k$ is an integer. Thus, there are finitely
many possibilities for $A_k$, so, passing to a subsequence, we
can assume that $A_k=A$ is a constant. Now we have by
(\ref{extra}),
\begin{align*}
y_k &= 2r n_k \th^{-N_k} = 2r\bigl(a_0 L_{N_k}^{(k)} +  \cdots +
a_{m-1} L_{N_K+m-1}^{(k)}\bigr)\th^{-N_k} + 2rA \th^{-N_k} \\
&= 2r z_0 (a_0 + a_1\th + \cdots + a_{m-1}\th^{m-1}) \\
&+ 2r \sum_{i=1}^M z_i \bigl(a_0 \th^{-\jk_i}+\cdots +
a_{m-1}\th^{-\jk_i+m-1}\bigr) + 2rA\th^{-N_k} \\
&+O\bigl(\th^{-N_k} \rho ^{N_k-\jk_M}\bigr),
\end{align*}
and finally, by (\ref{invert}),
\begin{equation}\label{equ3}
y_k=  z_0 + \sum_{i=1}^M z_i \th^{-\jk_i}  + 2rA\th^{-N_k} +
O\bigl(\th^{-N_k} \rho ^{N_k-\jk_M}\bigr).
\end{equation}

\begin{lemma} \label{lem-converge}
\begin{align*}
|\muhat(n_k)| &= \prod_{j=-\infty}^{N_k} |\cos(\pi y_k \Th^j)|\\
&\to \prod_{i=0}^M \prod_{j=-\infty}^\infty |\cos(\pi z_i
\Th^j)|\cdot \prod_{j=0}^\infty |\cos(2\pi rA \Th^{-j})|,\quad
k\to+\infty.
\end{align*}
\end{lemma}

Recall that by the construction of $z_i$ we have for $0 \le i \le
M$:
\begin{equation} \label{eq-dioph}
\|z_i \Th^j\| = O(\rho^j),\quad j\to \infty.
\end{equation}
It follows that the bi-infinite products in the right-hand side
converge. This lemma will clearly imply our theorem, as $z_i \in
\Q(\Th)$ and $A\in \Q$, so there are countable many possible
limits.

\medskip

\noindent {\em Proof of Lemma~\ref{lem-converge}.} We can find
integers $E_i^{(k)},\ i=0,\ldots,M$, so that
\begin{equation*}
\jk_i < E_i^{(k)} < \jk_{i+1},\quad\lim_{k\to\infty}
\bigl(E_i^{(k)}-\jk_i\bigr)
=\lim_{k\to\infty}\bigl(\jk_{i+1}-E_i^{(k)}\bigr)=+\infty.
\end{equation*}
We are going to show that
\begin{align} \label{conv1}
F_0(k)&:= \frac{\prod_{j=-\infty}^{E_0^{(k)}} |\cos(\pi y_k
\theta^j)|} {\prod_{j=-\infty}^\infty |\cos(\pi z_0 \Th^j)|} \to
1,\quad k\to \infty; \\
\label{conv2}F_i(k)&:= \frac{\prod_{j=E_{i-1}^{(k)}+1}^{E_i^{(k)}}
|\cos(\pi y_k \theta^j)|}{\prod_{j=-\infty}^\infty |\cos(\pi z_i
\Th^j)|} \to
1,\ \ \ k\to \infty, \ i=1,\ldots,M; \\
\label{conv3}F_{M+1}(k)&:= \frac{\prod_{j=E_M^{(k)}+1}^{N_k}
|\cos(\pi y_k \theta^j)|}{\prod_{j=0}^\infty |\cos(2\pi rA
\Th^{-j})|} \to 1,\quad k\to \infty.
\end{align}
These statements will imply the lemma.


First we verify (\ref{conv1}). Observe that
$\prod_{j=E_0^{(k)}}^\infty |\cos(\pi z_0 \Th^j)| \to 1$, since
$E_0^{(k)}\to \infty$ and the denominator in $F_0(k)$ converges.
Thus, it remains to show that
$$
\prod_{j=-\infty}^{E_0^{(k)}} \frac{|\cos(\pi y_k
\theta^j)|}{|\cos(\pi z_0 \theta^j)|} \to 1.
$$
\begin{sloppypar}
Note that $y_k = z_0 + O\bigl(\Th^{-\jk_1}\bigr)$ by (\ref{equ3}).
By assumption (\ref{eq-A}), $|\cos(\pi y_k\Th^j)| \ge \eta$ for $j
\le N_k$. Since $|y_k\Th^j - z_0 \Th^j|  =
O\bigl(\Th^{j-\jk_1}\bigr)$, we have $|\cos(\pi y' \Th^j)| \ge
\eta/2$ for $y'$ between $z_0$ and $y_k$, for all $j \le
E_0^{(k)}$, for $k$ sufficiently large. Then we can take logarithm
of each term and use the mean value theorem to get
\end{sloppypar}
\begin{align*}
\left|\log\left| \frac{\cos(\pi y_k \theta^j)}{\cos(\pi z_0
\theta^j)}\right| \right| &\le |\tan (\pi y'\theta^j)|\cdot |y_k
-z_0 | \Th^j \\
&\le \frac2\eta |y_k -z_0 | \Th^j = O\bigl(\Th^{j-\jk_1}\bigr),
\end{align*}
where $y'$ is between $z_0$ and $y_k$. Summing over $j \le
E_0^{(k)}$ and letting $k\to \infty$ yields the desired claim,
since $\jk_1-E_0^{(k)}\to+\infty$.


Now we verify (\ref{conv2}). Since $E_{i-1}^{(k)}-\jk_i\to
-\infty,\ E_i^{(k)}-\jk_i\to+\infty$, and the denominator in
$F_i(k)$ converges, it suffices to show that
\begin{equation}\label{eq-prod}
\prod_{j=E_{i-1}^{(k)}+1}^{E_i^{(k)}} \frac{|\cos(\pi y_k
\theta^j)}{\bigl|\cos\bigl(\pi z_i \theta^{j-\jk_i}\bigr)\bigr|}
\to 1.
\end{equation}
In view of (\ref{equ3}), we can write
$$
y_k \Th^j = z_0 \Th^j + \sum_{\ell=1}^{i-1} z_\ell
\Th^{j-\jk_\ell} + z_i \Th^{j-\jk_i} +
O\bigl(\Th^{j-\jk_{i+1}}\bigr),
$$
for $j = E_{i-1}^{(k)}+1,\ldots,E_i^{(k)}$. By (\ref{eq-dioph}),
for $j = E_{i-1}^{(k)}+1,\ldots,E_i^{(k)}$ we have
$$
y_k \Th^j =O\bigl(\rho^{j-\jk_{i-1}}\bigr) + z_i \Th^{j-\jk_i} +
O\bigl(\Th^{j-\jk_{i+1}}\bigr)\bmod\Z.
$$
For $k$ sufficiently large, the denominators in (\ref{eq-prod})
are bounded away from $0$, as above, when we checked
(\ref{conv1}), and we obtain
$$
\left|\log\left| \frac{\cos(\pi y_k \theta^j)}{\cos\bigl(\pi z_i
\theta^{j-\jk_i}\bigr)} \right|\right| \le \const\cdot
\bigl(\rho^{j-\jk_{i-1}} + \Th^{j-\jk_{i+1}}\bigr).
$$
Summing over $j = E_{i-1}^{(k)}+1,\ldots,E_i^{(k)}$, we obtain
that the logarithm of the product  in (\ref{eq-prod}) is bounded
in modulus by $\const\cdot\bigl(\rho^{E_{i-1}^{(k)} -\jk_{i-1}} +
\Th^{E_i^{(k)}-\jk_{i+1}}\bigr)$ which tends to 0, as $k\to
\infty$.

It remains to check (\ref{conv3}). Since $E_M^{(k)}-N_k \to
-\infty$ and the denominator in $F_{M+1}(k)$ converges, it is
sufficient to show that
\begin{equation} \label{eq-prod2}
\prod_{j=E_M^{(k)}+1}^{N_k} \frac{|\cos(\pi y_k
\theta^j)|}{|\cos(2\pi rA \theta^{j-N_k})|} \to 1.
\end{equation}
We have for $j = E_M^{(k)}+1,\ldots,N_k$, from (\ref{equ3}), in
view of (\ref{eq-dioph}):
$$
y_k \Th^j =2rA\Th^{j-N_k} + O\bigl(\rho^{j-\jk_M}\bigr) +
O\bigl(\Th^{j-N_k}\rho^{N_k-\jk_M}\bigr)\bmod\Z.
$$
For $k$ sufficiently large, the denominators in (\ref{eq-prod2})
are bounded away from $0$, as above, when we checked
(\ref{conv1}), and we can write
$$
\left|\log\left| \frac{\cos(\pi y_k \theta^j)}{\cos(2\pi
rA\theta^{j-N_k})} \right|\right| \le \const\cdot (\rho^{j-\jk_M}
+ \Th^{j-N_k}\rho^{N_k-\jk_M}).
$$
Summing over $j = E_M^{(k)}+1,\ldots,N_k$, we obtain that the
logarithm of the product  in (\ref{eq-prod2}) is bounded above in
modulus by $\const\cdot(\rho^{E_M^{(k)} -\jk_M} +
\rho^{N_k-\jk_M})$ which tends to 0, as $k\to \infty$. This
concludes the proof of Lemma~\ref{lem-converge}.
Inequality~(\ref{eq:loeq}) and thus, Theorem~\ref{thm2} and
Theorem~\ref{main} are proved as well. \qed

\section{Proofs of other results}
\label{sec:proofs}

\subsection{Proof of Theorem~\ref{thm2}~(ii)} Let $J_\th$ denote
the set of limit points of $\{\muhat_\th(t): t>0\}$ as
$t\to\infty$.
\begin{lemma}\label{lem:int} $J_\th$ is a non-empty interval.
\end{lemma}
\begin{proof}
By the theorem of Erd\H{o}s \cite{E}, $J_\th$ contains a
non-zero point. On the other hand, $0\in J_\th$, since
$\muhat_\th(\th^n/4)=0,\ n\ge1$. Let $a = \inf J_\th,\  b = \sup
J_\th$. Then there are sequences $u_i,v_i\to \infty$ such that
$\muhat_\th(u_i) \to a$ and $\muhat_\th(v_i) \to b$. Without
loss of generality, $u_i < v_i < u_{i+1}$ for all $i$. Since
$t\mapsto\wh{\mu_\th}(t)$ is continuous, for any $\e>0$, for all
$i$ sufficiently large, any value between $a+\e$ and $b-\e$ is
assumed by $\muhat_\th(t)$ at least once in $(u_i,v_i)$. Thus,
$J_\th = [a,b]$.
\end{proof}

Our goal is to prove that for a.e.\ $r>0$, $J_\th$ is in fact
the set of limit points for the sequence $\{\muhat_\th(rn): n\in
\Z\}$ as well.

Let $\{y_k\}_{k\ge 1}$ be a sequence dense in $J_\th$. By the
definition of $J_\th$, for any $k\ge 1$, there is a sequence
$t_i^{(k)}\to +\infty$ as $i\to\infty$, such that
$\lim_{i\to\infty}\muhat_\th\bigl(t_i^{(k)}\bigr)= y_k$. Recall
the following well-known fact.

\begin{prop} \label{prop-distr} \cite[Chap.~1, Sec.~4, Cor~4.3]{KN}
For any unbounded sequence $\{x_i\}_{i\ge 1}$, the set $\{\alpha
x_i\}_{i \ge 1}$ is dense modulo 1 for Lebesgue-a.e. $\alpha$.
\end{prop}

Thus, for a.e.\ $r>0$, the sequence
$\bigl\{r^{-1}t_i^{(k)}\bigr\}_{i\ge 1}$ is dense modulo~1 for
all $k\ge 1$. Fix such an $r$. Then for any $k\ge 1$, there is a
subsequence $t_{i_j}^{(k)}$ such that $r^{-1}t_{i_j}^{(k)} \to
0$ mod 1 (of course, $i_j$ may depend on $k$). Thus, for any
$k$, there exist $n_j \in \Nat$ such that $rn_j - t_{i_j}^{(k)}
\to 0$, as $j\to \infty$. Note that $\left| \frac{d}{dt}
\wh{\mu_\th}(t) \right| \le C$, since $\mu_\th$ has compact
support on $\R$, whence \be \label{eq-lima} \{\wh{\mu_\th}(a_n):
n \in \Nat\}' = \{\wh{\mu_\th}(b_n): n \in \Nat\}' \ee for any
$a_n, b_n \to \infty$, with $a_n - b_n \to 0$. Therefore,
$$
\lim_{j\to\infty} \wh{\mu_\th}(rn_j) = \lim_{j\to\infty}
\wh{\mu_\th}\bigl(t_{i_j}^{(k)}\bigr) = y_k.
$$
It follows that for a.e. $r>0$, the set of limit points of
$\{\muhat_\th(rn): n\in \Nat\}$ contains all $y_k$, which are
dense in $J_\th$, and hence it contains all of $J_\th$. \qed

\subsection{Proof of Theorem~\ref{thm:deriv}} In view of
Theorem~\ref{thm2}~(i) and the fact that
$\Fk_{\mu_\th,r}^{(n+1)}\subset\Fk_{\mu_\th,r}^{(n)},\ n\ge1$, it
suffices to show that the $n$'th derived set of limit points of
$\Fk_{\mu_\th,r}$ is {\em at least} countable for $r\in \Q(\th)$.

Since $r \in \Q(\th)$, there exists $p \in \Z[x]$ such that $\La=r
p(\th) \in \Z[\th]$. Let $q\in\Z[\th]$ be an arbitrary number. We
have
\begin{equation*}
\lim_{k\to\infty} \wh{\mu_\th}(r \langle p(\th) q\th^k\rangle ) =
\lim_{k\to\infty} \wh{\mu_\th}(\La q\th^k)=
\prod_{j=-\infty}^\infty \cos(2\pi \La q\th^j).
\end{equation*}
The first equality holds by (\ref{eq-lima}), as $\|h\th^n\| =
O(\rho^n)=o(1)$ for any $h\in \Z[\th]$. The second equality
follows from (\ref{eq-four}). Put
\[
\phi_\La(q)=\prod_{j=-\infty}^\infty \cos(2\pi \La q \th^j)
\]
and
\[
\Om_\La:=\left\{\phi_\La(q) : q\in\Z[\th]\right\}.
\]
We just proved that $\Om_\La \subset \Fk_{\mu_\th,r}'$, so
$\Om_\La^{(n)}\subset\Fk_{\mu_\th,r}^{(n+1)},\ n\ge1$. Our next
goal is to show first that $\Om_\La'$ (and hence
$\Fk_{\mu_\th,r}''$) is infinite for every $\La\in\Z[\th]$. Let
$a,b\in\Z[\th]$ and put $q_n(a,b):=a+b\th^n$. Then
\begin{align*}
\phi_\La(q_n(a,b))&=\prod_{j=-\infty}^{\left[\frac n2\right]}
\cos(2\pi \La(a+b\th^n)\th^{-j})\cdot
\prod_{\left[\frac{n}2\right]}^{{\left[\frac{3n}2\right]}}
\cos(2\pi \La(a+b\th^n) \th^{-j})\\
&\cdot \prod_{j=\left[\frac{3n}2\right]+1}^\infty\cos(2\pi \La
(a+b\th^n) \th^{-j}),
\end{align*}
and similarly to Lemma~\ref{lem-converge}, it is easy to see
that the first product tends to $\phi_\La(a)$, the second one
tends to $\phi_\La(b)$ and finally, the last one tends to 1.
Hence
\begin{equation}
\phi_\La(q_n(a,b))\to\phi_\La(a)\phi_\La(b),\quad
n\to+\infty.\label{eq:produ}
\end{equation}
Recall that $\th\neq2$ (when $\phi_\La(q)\equiv0$ for any $\La$
and $q$), so there always exists $q$ such that
$0<|\phi_\La(q)|<1$, whence $\Om_\La$ is infinite. Furthermore,
since $q_n(a,b) \in \Z[\th]$, (\ref{eq:produ}) implies that
$\Om_\La^2\subset\Om_\La'$, where $\Om_\La^2=\{\om_1\om_2 :
\om_i\in\Om_\La\}$. Therefore, $\Om_\La'$ is infinite as well.
Now, $\Om_\La^4\subset(\Om_\La')^2 \subset
(\Om_\La^2)'\subset\Om_\La''$, whence $\Om_\La''$ is also
infinite, etc. By induction, $\Fk_{\mu_\th,r}^{(n)}$ is
countable for each $n\ge1$. \qed

\section{Concluding remarks}

\noindent \textbf{1.} Our first remark concerns the proof of
Theorem~\ref{thm2}~(ii) (see the beginning of the previous
section). In fact, what we use is the following

\begin{lemma} Assume $f\in C(\R_+)\cap L^\infty(\R_+)$, and $J$
is the set of the limit points of $f$ as $t\to+\infty$. Then the
derived set for $\{f(rn) : n\in\Nat\}$ as $n\to+\infty$ is equal
to $J$ for a.e. $r>0$.
\end{lemma}

The proof of this lemma is exactly the same as above for
$f(t):=\wh{\mu_\th}(t)$. This claim is probably known but
we did not find it in the literature.

\medskip\noindent
\textbf{2.} Our second remark consists in a simple observation
that the expression for the limit points of $\Fk_{\mu_\th,r}$ in
Lemma~\ref{lem-converge} (without the moduli) is in fact a
general formula for $x\in\Fk_{\mu_\th,r}'$. More precisely, let
\[
\P_\th=\{\xi : \|\xi\th^n\|\to0,\ n\to+\infty\}.
\]
As is well known, $\Z[\th]\subsetneq\P_\th\subsetneq\Q(\th)$
(see, e.g., \cite{Cas}), and $\P_\th$ is obviously a group under
addition. Then our claim is that for any $M\in\Z_+,
(z_i)_{i=0}^M\in\P_\th^{M+1}$ and $A\in\Z$,
\[
x=\prod_{i=0}^M \prod_{j=-\infty}^\infty \cos(\pi z_i
\th^j)\cdot \prod_{j=0}^\infty \cos(2\pi rA
\th^{-j})\in\Fk_{\mu_\th,r}'.
\]
Indeed, put
\[
n_k=\langle
(2r)^{-1}(z_0\th^{(M+1)k}+z_1\th^{Mk}+\dots+z_M\th^k)\rangle+A.
\]
The proof is practically the same as for
Lemma~\ref{lem-converge}, and we leave it to the reader.

\medskip\noindent \textbf{3.} Our next remark concerns translations
of Bernoulli convolutions. Let $\ga\in\R$; then shifting the
origin by $\ga$ results in multiplying $\wh{\mu_\th}(rn)$ by
$e^{2\pi i \ga n}$.

\begin{prop} The closure of the set 
$\{\wh{\mu_\th}(nr)e^{2\pi i \ga n} :\, n\in\Z\}$ is countable for
$r\in\Q(\th)$ and $\ga\in\Q(\th)$, but uncountable for $r\in\Q(\th)$ and
a.e.\ $\ga\in \R$.
\end{prop}
\begin{proof} First suppose that $r\in\Q(\th)$ and $\ga\in\Q(\th)$.
It follows from the proof of the main theorem that
if $|\muhat_\th(rn_k)| \not\to 0$, then the formula~(\ref{equ3}) gives
a general expression for $n_k$, with $k$ sufficiently large.
Now it is enough to note that for any $\xi \in \Q(\th)$, the sequence
$\{\|\xi\th^n\| : n\in\Nat\}$ has finitely many limit points\footnote{This
was known to Pisot, see e.g.\ \cite[p.\,96]{Pisot}} and hence
the sequence $\{\|\ga n_k\|:\ k\in\Nat\}$ has finitely many limit points.

On the other hand, given $r\in\Q(\th)$, we can fix $n_k$ so that
$\muhat_\th(rn_k)\to a \ne 0$. Using Proposition~\ref{prop-distr} again,
we see that for a.e.\ $\ga\in \R$, the sequence
$\{\ga n_k : k\in\Nat\}$ is dense modulo 1, and therefore,
the set of limit points of  the sequence $\{\muhat_\th(rn_k)e^{2\pi i \ga n_k}
\}$ is the circle of radius $|a|$.
\end{proof}



\medskip\noindent {\bf 4.} Denote by $\mu_\th(r,\ga)$ the measure on the
circle whose Fourier coefficients are $\wh{\mu}_\th(rn)e^{2\pi in\ga}$
(that is, the translation of a scaled copy of $\mu_\th$).
We have shown that for ``most'' $(r,\ga)$ the spectrum
$\mbox{sp}(T_{\mu_\th(r)},L^2)$ contains a continuum
(an interval or a circle). In these cases we cannot use Sarnak's result
\cite{Sar} to claim that the spectra are the same in all $L^p$, for
$p \in (1,\infty)$. However, the remarks in \cite{Christ} indicate
that $\mu_\th(r,0)$ is $L^p$-improving for any $r>1$ and $\th>1$,
and any translation of the measure, obviously, preserves the property.
Thus, by \cite[Th.~4.1]{GHR} we see that the claim on coincidence of
spectra in all $L^p$ is valid for all $\mu_\th(r,\ga)$.

%
%

\subsection*{Acknowledgment} We would like to thank
the Erwin Schr\"odin\-ger Institute for Mathematical Physics at
Vienna as well as Peter Grabner and Robert Tichy, the organizers
of the workshop ``Arithmetics, Automata and Asymptotics" in Graz
(July, 2002) for their hospitality.

\end{document}